# MISSING AT RANDOM, LIKELIHOOD IGNORABILITY AND MODEL COMPLETENESS[1]

By Guobing Lu and John B. Copas

*University of Bristol and University of Warwick*

This paper provides further insight into the key concept of missing at random (MAR) in incomplete data analysis. Following the usual selection modelling approach we envisage two models with separable parameters: a model for the response of interest and a model for the missing data mechanism (MDM). If the response model is given by a complete density family, then frequentist inference from the likelihood function ignoring the MDM is valid if and only if the MDM is MAR. This necessary and sufficient condition also holds more generally for models for coarse data, such as censoring. Examples are given to show the necessity of the completeness of the underlying model for this equivalence to hold.

**1. Introduction.** A full parametric model for missing data comprises two components: one is for the complete data and the other is for the *missing data mechanism* (MDM). The former describes the probability distribution that governs the data generation process of interest, while the latter characterizes the observation process by which some data may be missing. The parameterizations of these two processes are often assumed to be separable, and our target is to make inference about the parameters involved in the complete data model using only the available incomplete data.

In practice, modelling incomplete data is a very difficult task since in most cases the incomplete data themselves contain little or no information about the MDM. The fundamental and most widely used assumption about the MDM is that it is a *missing at random* (MAR) model [Rubin (1976)]. The basic idea is that the probability that a response variable is observed can depend only on the values of those other variables which have been

Received April 2002; revised May 2003.

[1]Supported by UK Engineering and Physical Sciences Research Council (EPSRC) Grant GR/L93980.

observed. This concept has been extensively studied, and effective computational methods for handling missing data under the MAR assumption have been well developed, for example, using the EM algorithm. Good references include Tanner (1993), Schafer (1997), Kenward and Molenberghs (1998) and Little and Rubin (2002) among many others.

A closely related, but logically distinct, concept is *ignorability*. The basic idea here is that inference based on the joint specification of both complete data and MDM models is the same as the inference we would obtain if we used the complete data model only, simply integrating out the values of any variables which are missing. It is well known that MAR, together with the assumption of separable parameters, is a sufficient condition for ignorability of the MDM in likelihood based inference. It is, however, not a necessary condition.

Although these concepts have been widely discussed, there have been some inconsistencies between different authors on how they are defined and interpreted, and in the choice of terminology. The Weblist impute@utdallas.edu gives an interesting summary of views. We avoid ambiguities by giving some precise definitions in Section 2.

In Section 3 we show that for models given by a complete family of distributions, MAR is both necessary and sufficient for ignorability. The result depends on a heritable property of completeness: that, with suitable reparameterizations, completeness of a multivariate distribution implies completeness of all conditional and marginal distributions. Examples are given to show that, for inference in a family of distributions which is not complete, an MDM can be ignorable without being MAR.

This necessary and sufficient condition is extended in Section 4 to the wider concept of *coarsening at random* introduced by Heitjan and Rubin (1991). Here, ideas for missing data are generalized to other kinds of incomplete data such as censoring or rounding.

Section 5 offers some concluding remarks.

## 2. Missing data and likelihood ignorability.

Let $Y = (Y_1, \ldots, Y_k)^T$ be a $k$-dimensional random vector with probability density function $f(y; \theta)$ on $\mathcal{Y} \subset \mathbb{R}^k$, where $\theta \in \Theta$ is a $d$-dimensional parameter of interest. Suppose that the observation process of $Y$ suffers from missing data and hence, associated with $Y$ there is also a binary random vector $R = (R_1, \ldots, R_k)^T$ indicating the observational status of $Y$, where $R_i$ takes the value 0 when the observation of $Y_i$ is missing and the value 1 when $Y_i$ is observed, $i = 1, \ldots, k$. Denote the range of $R$ by

$$\mathcal{R} = \{(r_1, \ldots, r_k) : r_i = 0 \text{ or } 1, \ i = 1, \ldots, k\} = \{0, 1\}^k.$$

We assume that the parameterization of the joint distribution of $Y$ and $R$ can be put into the selection model form

$$(2.1) \qquad f(y, r; \theta, \psi) = f(y; \theta) f(r | y, \psi), \qquad (\theta, \psi) \in \Theta \times \Psi,$$

in which the parameters $\theta$ and $\psi$ are assumed to be distinct [Rubin (1976)]. The conditional density $f(r|y,\psi)$ characterizes the probabilistic relation between the data-observation process and the values of the data themselves and hence specifies a model for the MDM. The joint distribution of $Y$ and $R$ can also be written in the pattern mixture form [Little (1994)] in which we model instead the conditional distribution of $Y$ given $R$, but the parameterization in (2.1) makes the MAR condition more transparent for the present discussion.

The pair of random variables $(Y, R)$ induces an observable random variable $Z$. Using a notation analogous to that for the coarsening function as defined in Heitjan and Rubin (1991), $Z$ is

$$Z = Z(Y, R) = (Z_1, \ldots, Z_k)^T,$$

$$\text{(2.2)} \qquad \text{where } Z_i = \begin{cases} Y_i, & \text{if } R_i = 1, \\ \mathbb{R}, & \text{if } R_i = 0, \end{cases} \quad i = 1, \ldots, k.$$

For notational convenience we allow the symbol $\mathbb{R}$ to appear in any position in the vector argument of a multivariate density function, using it to denote the marginal density of the other variables. For example, suppose $f(t_1, t_2)$ is a density on $\mathbb{T}_1 \times \mathbb{T}_2 \subset \mathbb{R}^2$ and $f_i(t_i)$, $i = 1, 2$, are the marginal densities. Then we identify $f(t_1, \mathbb{R})$ with $f_1(t_1)$ and $f(\mathbb{R}, t_2)$ with $f_2(t_2)$. Trivially, $f(\mathbb{R}, \mathbb{R}) \equiv 1$. With this convention the density of $Z$ can be expressed as

$$f(z; \theta, \psi) = \int f(y; \theta) f(r|y; \psi) \, dy^{(\mathbb{1}-r)}$$

$$\text{(2.3)}$$

$$= f(y^{(r)}; \theta) \int f(y^{(\mathbb{1}-r)}|y^{(r)}; \theta) f(r|y; \psi) \, dy^{(\mathbb{1}-r)},$$

where $\mathbb{1}$ is the $k$-dimensional vector with all elements equal to 1, $y^{(r)}$ and $y^{(\mathbb{1}-r)}$ are, respectively, the observed subvector and the missing subvector of $y$ given by

$$y^{(r)} = (y_i : r_i = 1, \ i = 1, \ldots, k)^T$$

and

$$y^{(\mathbb{1}-r)} = (y_i : r_i = 0, \ i = 1, \ldots, k)^T,$$

and for each variable $y_i$ contained in $y^{(\mathbb{1}-r)}$ the integral in (2.3) is over its whole range.

In the above setting Rubin's MAR condition [Rubin (1976)] can be expressed as follows.

DEFINITION 2.1. A MDM is said to be MAR if the conditional distribution $f(r|y; \psi)$ has the special form

$$\text{(2.4)} \qquad f(r|y; \psi) = h_r(z(y, r); \psi) \qquad \text{for all } (y, r) \in \mathcal{Y} \times \mathcal{R},$$

where, for any fixed $\psi$ and $r$, $h_r(\cdot; \psi)$ is a function mapping $\mathbb{R}^{(\mathbb{1}^T r)}$ into $[0, 1]$.

Note that the dimension of the space $\mathbb{R}^{(\mathbb{1}^T r)}$ varies with the value of $r$, and hence $h_r(\cdot; \psi)$ is a family of $2^k$ functions indexed by the subscript $r$. Under MAR the MDM depends on $(y, r)$ only through the function $z(y, r)$, that is, through the observed part of the sample $y$.

A MDM model is ignorable for inference about the parameter $\theta$ if the inference based on the combination of both the complete data model and the MDM model coincide with the inference based on the complete data model alone. For likelihood based inference we assume that ignorability is an intrinsic property of the joint model $f(y, r; \theta, \psi)$ rather than a property of any specific sample realization. Thus we are interested in frequentist inference from the likelihood function, rather than inference from the particular likelihood function we get from the observed sample. To emphasize this we use the term *likelihood ignorable* (LIG) in the following definition.

DEFINITION 2.2.   A MDM is said to be LIG if the integral

$$(2.5) \qquad \int f(y^{(\mathbb{1}-r)} | y^{(r)}; \theta) f(r | y; \psi) \, dy^{(\mathbb{1}-r)}$$

is free of $\theta$ for almost all realizations of $(y, r) \in \mathcal{Y} \times \mathcal{R}$ and for all $(\theta, \psi) \in \Theta \times \Psi$.

The contribution of observation $z$ to the likelihood is the product of the two terms in the right-hand side of (2.3). LIG means that the second term [the integral over $y^{(\mathbb{1}-r)}$] does not affect the likelihood as far as inference about $\theta$ is concerned. Equivalently, the contribution of this second term of the log likelihood disappears when we differentiate with respect to $\theta$. All that matters is the first term, which is just the marginal joint density of those components of $Y$ which are actually observed.

Notice that MAR is a property of the conditional distribution $f(r | y; \psi)$, whereas LIG depends on both $f(r | y; \psi)$ *and* the response model $f(y; \theta)$.

Under the MAR model,

$$\int f(y^{(\mathbb{1}-r)} | y^{(r)}; \theta) f(r | y; \psi) \, dy^{(\mathbb{1}-r)}$$

$$= \int f(y^{(\mathbb{1}-r)} | y^{(r)}; \theta) h_r(z(y, r); \psi) \, dy^{(\mathbb{1}-r)}$$

$$= h_r(z(y, r); \psi),$$

which is independent of $\theta$. Hence MAR is a sufficient condition for LIG. We seek the conditions under which MAR is also a necessary condition for LIG.

**3. Necessary and sufficient conditions for LIG.** In this section we show that if the family of density functions $f(y; \theta)$ forms a complete class, then MAR is both necessary and sufficient for LIG.

First some preliminaries about completeness. Recalling the elementary definition [e.g., Zacks (1971)], a family of probability density functions $\{f(y; \theta) : \theta \in \Theta\}$ on $\mathcal{Y} \subset \mathbb{R}^k$ is said to be complete if the identity

$$\int t(y) f(y; \theta) \, dy = 0 \qquad \text{for all } \theta \in \Theta$$

implies that $t(y) = 0$ for almost all $y \in \mathcal{Y}$.

Now let $Y^{(1)}$ be a subvector of the random vector $Y$, and let $Y^{(2)}$ be the corresponding complementary subvector. Denote the sample spaces of $Y^{(i)}$ by $\mathcal{Y}^{(i)}$, $i = 1, 2$. Then the joint family $f(y; \theta)$ can be decomposed into

$$f(y; \theta) = f_1(y^{(1)}; \theta) f_{2|1}\{y^{(2)} | y^{(1)}; \theta(y^{(1)})\},$$

where $\theta(y^{(1)})$ is a function: $\mathcal{Y}^{(1)} \mapsto \Theta$ [see Arnold, Castillo and Sarabia (1999)]. We remark that, in general, even if the joint density family can be identified by the parameter $\theta$, neither the marginal density family $f_1(y^{(1)}; \theta)$ nor the conditional family $f_{2|1}\{y^{(2)} | y^{(1)}; \theta(y^{(1)})\}$ is identified by the same parameter. However, there is always a many-to-one function: $\phi_1 : \Theta \mapsto \Phi_1 \subset \Theta$ such that

$$\{f_1(y^{(1)}; \phi_1) : \phi_1 \in \Phi_1\} = \{f_1(y^{(1)}; \theta) : \theta \in \Theta\},$$

and the new parameter $\phi_1(\theta)$ is identified. Similarly, for any given $y^{(1)}$, the conditional family $f_{2|1}$ can be identified by $\phi_2(\theta; y^{(1)})$. Detailed discussion of the problems of reparameterization and identification will in general call for a topological group structure in the parameter space $\Theta$, but for the purpose of describing completeness we merely borrow the form of the parameterization to make the representation clear.

The following lemma says that completeness is a heritable property from the joint density family to its marginal and conditional density families.

LEMMA 3.1. *Suppose that* $\{f(y; \theta) : \theta \in \Theta\}$ *is a complete density family. Then the following hold:*

(a) *the marginal family* $[f\{y^{(1)}; \phi_1(\theta)\} : \theta \in \Theta]$ *is complete;*
(b) *for almost all* $y^{(1)} \in \mathcal{Y}^{(1)}$ *the conditional families* $[f_{2|1}\{y^{(2)} | y^{(1)}, \phi_2(\theta; y^{(1)})\} : \theta \in \Theta]$ *are complete.*

See the Appendix for the proof of Lemma 3.1.

Now we apply Lemma 3.1 to the conditional density family $f(y^{(1-r)} | y^{(r)}; \theta)$ to yield the following theorem.

THEOREM 3.1. *For the selection model* (2.1) *assume that* $\{f(y; \theta) : \theta \in \Theta\}$ *is a complete family. Then the necessary and sufficient condition for LIG is that the MDM is MAR.*

PROOF. We only need to verify that LIG implies MAR. From Definition 2.2 LIG implies that the integral (2.5) is independent of $\theta$ for any given $r$ and almost all $y^{(r)}$. Denoting its values by $w(y^{(r)}, r; \psi)$, we have the equality

$$(3.1) \qquad \int \{f(r|y; \psi) - w(y^{(r)}, r; \psi)\} f(y^{(\mathbb{1}-r)}|y^{(r)}; \theta) \, dy^{(\mathbb{1}-r)} = 0.$$

Because of the inheritance property of completeness, $f(y^{(\mathbb{1}-r)}|y^{(r)}; \theta)$ are complete density families for all $r \in \mathcal{R}$, and hence, for any values of $r$ and $\psi$, (3.1) implies

$$f(r|y; \psi) = w(y^{(r)}, r; \psi).$$

Thus the MDM $f(r|y; \psi)$ depends on $y$ only through $y^{(r)}$ and so must have the form of $h_r(z(y, r); \psi)$ in (2.4). That is, the MAR condition holds. □

The following examples show that, for an incomplete density family, LIG does not guarantee MAR.

EXAMPLE 3.1. Consider the bivariate normal density family

$$Y = \begin{pmatrix} Y_1 \\ Y_2 \end{pmatrix} \sim N \left\{ \begin{pmatrix} \theta \\ \theta/2 \end{pmatrix}, \begin{pmatrix} 1 & 1/2 \\ 1/2 & 1 \end{pmatrix} \right\}, \qquad \theta \in \mathbb{R}.$$

Clearly this is not a complete family since $E(Y_1 - 2Y_2) = 0$ for all values of $\theta$.

Suppose that $Y_1$ is always observed but $Y_2$ may be missing. The MDM is then characterized by the functions

$$h_{(1,0)}(y; \psi), \qquad h_{(1,1)}(y; \psi) = 1 - h_{(1,0)}(y; \psi), \qquad h_{(0,0)} = h_{(0,1)} = 0.$$

The MAR condition demands that $h_{(1,1)}(y; \psi)$ as a function of $y = (y_1, y_2)$ is independent of $y_2$ for all $\psi \in \Psi$. However, in this example the conditional density of $Y_2$ given $Y_1$ is independent of $\theta$. Hence, for $r = (1, 0)$ and any arbitrary function $h_{(1,0)}(y; \psi)$, the integral

$$\int f(y^{(\mathbb{1}-r)}|y^{(r)}; \theta) f(r|y; \psi) \, dy^{(\mathbb{1}-r)} = \int f(y_2|y_1) h_{(1,0)}(y; \psi) \, dy_2$$

does not depend on $\theta$. Thus in this case *any* MDM is LIG.

Example 3.2. Now extend Example 3.1 by supposing that $Y_1$ or $Y_2$, but not both, may be missing. Suppose that the MDM is

$$f(r|y; \psi) = \begin{cases} h_{(1,1)}(y; \psi), & \text{when } r = (1,1), \\ h_{(1,0)}(y; \psi), & \text{when } r = (1,0), \\ h_{(0,1)}(y_2; \psi), & \text{when } r = (0,1), \\ h_{(0,0)} \equiv 0, & \text{when } r = (0,0), \end{cases}$$

where

$$\sum_{i,j} h_{(i,j)}(y; \psi) = 1$$

for all $y$ and $\psi$. Since $h_{(1,0)}$ depends on both $y_1$ and $y_2$, the MDM is not MAR. However, because $h_{(1,1)}, h_{(0,1)}$ and $h_{(0,0)}$ satisfy the MAR condition, we only need to check the LIG condition for $r = (1,0)$. However, this is just the same as Example 3.1, so the LIG condition holds.

Example 3.3. Let $Y = (Y_1, Y_2, \ldots, Y_k)^T$ be i.i.d. $N(\theta, 1)$. Then $S = \frac{1}{k} \sum_{i=1}^{k} Y_i$ is a sufficient statistic and the vector of sample differences

$$A = (Y_1 - Y_2, Y_2 - Y_3, \ldots, Y_{k-1} - Y_k)^T$$

is an ancillary statistic for $\theta$. These statistics are independent, so we have

$$f(y; \theta) = f(y|s) f(s; \theta) = f(a) f(s; \theta).$$

Similarly, for any given $r$ with $\mathbb{1}^T r < k - 1$, we can define the corresponding statistics $s_r$ and $a_r$ for the subvector $y^{(\mathbb{1}-r)}$.

Now suppose that the MDM takes the form

$$f(r|y; \psi) = h_r(a_r, y^{(r)}; \psi).$$

Clearly this is not MAR because $h_r$ depends on $y^{(\mathbb{1}-r)}$ through $a_r$. However,

$$\int f(y^{(\mathbb{1}-r)}|y^{(r)}; \theta) f(r|y; \psi) \, dy^{(\mathbb{1}-r)}$$

$$= \int\int f(a_r) f(s_r; \theta) h_r(a_r, y^{(r)}; \psi) \, ds_r \, da_r$$

$$= \int f(a_r) h_r(a_r, y^{(r)}; \psi) \, da_r,$$

which does not depend on $\theta$. Hence this MDM is LIG.

**4. Extension to coarsening at random.**  The coarse data model of Heitjan and Rubin (1991) is a more general way of describing incomplete data problems. Here $Z$, the observable outcome, is a measurable subset of the sample space, such as a half line (when a life time is known to exceed a censoring time) or a finite interval (when an observation is rounded). The notion of *coarsening at random* (CAR) was introduced by Heitjan and Rubin (1991) as a natural extension of MAR to coarse data and was further studied by Heitjan (1993, 1994, 1997) and Jacobsen and Keiding (1995).

Following Heitjan and Rubin (1991), a random variable $G$, the so-called coarsening variable, defines the measurable subset $Z$ as $Z = Z(Y, G)$. Equation (2.2) is the special case of this when $G = R$. The conditional distribution of $G$ given $Y$ defines the coarsening data mechanism (CDM)

$$f(g|y; \psi) = h(g, y; \psi),$$

where, again, the parameter $\psi$ is assumed to be distinct from the parameter $\theta$ in the main model $f(y; \theta)$.

With the CDM $h(g, y; \psi)$, the conditional distribution of $Z$ given $y$ can be expressed as

$$\kappa(z, y; \psi) = \int_{\{g \,:\, Z(y,g)=z\}} h(g, y; \psi) \, dg.$$

For a rigorous expression for this conditional density in the case when $G$ is continuous, see Jacobsen and Keiding (1995).

The following definition of CAR is due to Heitjan and Rubin (1991).

DEFINITION 4.1.   The CDM is CAR if, for any fixed observed subset $z$, and for each value of $\psi$, $\kappa(z, y; \psi)$ takes the same value for all $y \in z$.

The likelihood function for $\theta$ based on an observed $z$ is proportional to the probability that $(Y, G)$ falls in the set $\{(y, g) : Z(y, g) = z\}$, which can be written as

$$\int_z \int_{\{g \,:\, Z(y,g)=z\}} f(y; \theta) h(g, y; \psi) \, dg \, dy = \int_z f(y; \theta) \kappa(z, y; \psi) \, dy.$$

This leads to the following definition.

DEFINITION 4.2.   The CDM is said to be LIG if, as functions of $\theta$,

$$(4.1) \qquad \int_z f(y; \theta) \kappa(z, y; \psi) \, dy \propto \int_z f(y; \theta) \, dy.$$

The generalization of Theorem 3.1 is as follows.

THEOREM 4.1.   *If $\{f(y; \theta) : \theta \in \Theta\}$ is a complete family, then a necessary and sufficient condition for LIG is that the CDM is CAR.*

That CAR implies LIG is immediate. For the converse, if the CDM satisfies (4.1), there exists $w(z; \psi)$ such that

$$\int_z f(y; \theta)\{\kappa(z, y; \psi) - w(z; \psi)\}\, dy = 0.$$

We now need a straightforward extension of Lemma 3.1, that if $f(y; \theta)$ is complete, then so is the conditional distribution of $y$ given $y \in z$ (the proof follows lines similar to the proof of Lemma 3.1 in the Appendix). This implies that $\kappa(z, y; \psi) = w(z; \psi)$ and hence the CDM is CAR.

The following example shows the necessity of model completeness for the equivalence of CAR and LIG when the coarsening variable $G$ is a continuous random variable.

EXAMPLE 4.1. Let $Y = (Y_1, Y_2)^T$ be the logarithms of two life times, assumed to follow the (incomplete) bivariate normal distribution in Example 3.1. Suppose that $Y_1$ is always observed but $Y_2$ suffers from censoring, with $G$ the corresponding (log) censoring time in a competing risks framework. The coarsening function is

$$Z(y, g) = \begin{cases} \{y\}, & \text{if } g \geq y_2, \\ \{y_1\} \times (g, \infty), & \text{if } g < y_2. \end{cases}$$

Suppose further that $(Y, G)$ are jointly Gaussian

$$\begin{pmatrix} Y_1 \\ Y_2 \\ G \end{pmatrix} \sim N\left\{ \begin{pmatrix} \theta \\ \theta/2 \\ \psi \end{pmatrix}, \begin{pmatrix} 1 & 1/2 & 0 \\ 1/2 & 1 & 1/2 \\ 0 & 1/2 & 1 \end{pmatrix} \right\}, \qquad \theta, \psi \in \mathbb{R}.$$

For this model we find

$$h(g, y; \psi) = \phi\left( \frac{g + (1/3)y_1 - (2/3)y_2 - \psi}{\sqrt{2/3}} \right)$$

and

$$\kappa(z, y; \psi) = \begin{cases} \int_{(y_2, \infty)} \phi\left( \dfrac{g + (1/3)y_1 - (2/3)y_2 - \psi}{\sqrt{2/3}} \right) dg, \\ \hspace{7cm} \text{if } z = \{y\}, \\ \phi\left( \dfrac{g + (1/3)y_1 - (2/3)y_2 - \psi}{\sqrt{2/3}} \right), \quad \text{if } z = \{y_1\} \times (g, \infty), \end{cases}$$

where $\phi(\cdot)$ is the standard normal density function. Clearly, $\kappa(z, y; \psi)$ does not take the same value for all $y \in z$ for each value of $\psi$, and so the CDM is not CAR. However, it is LIG, because for an observation $z^o = \{y_1^o\} \times (g^o, \infty)$,

$$\int_{z^o} f(y_1^o, y_2; \theta)\kappa(z^o, y; \psi)\, dy$$

$$= \int_{(g^o,\infty)} f(y_1^o;\theta)f(y_2|y_1^o)\phi\left(\frac{g^o+(1/3)y_1^o-(2/3)y_2-\psi}{\sqrt{2/3}}\right)dy_2$$

$$\propto \int_{(g^o,\infty)} f(y_1^o,y_2;\theta)\,dy_2,$$

as the conditional density of $y_2$ given $y_1$ is independent of $\theta$.

## 5. Remarks.

REMARK 5.1.   In missing data analysis we will usually assume that the data arise from $n$ i.i.d. realizations from the joint distribution of $(Y,R)$ in (2.1). If the distribution of $Y$ is complete, asking whether the MDM affects inference (in the sense of LIG) is then equivalent to asking whether $R$ depends on unobserved components of $Y$. In general, however, MAR is a stronger requirement than LIG. In Example 3.3, for instance, $Y$ itself takes the form of an i.i.d. sample, but the components of $R$ may be dependent. Now $R$ can depend in a arbitrary way on ancillary statistics without upsetting inference about $\theta$.

REMARK 5.2.   Many familiar statistical models used in practice involve replication and i.i.d. residuals and are not complete, such as Example 3.3. In normal linear models more generally, ignorable MDMs can still depend on the standardized sample residuals.

REMARK 5.3.   When covariates, say $X$, are involved in incomplete data analysis, we may wish to model conditionally on $X$, and hence the equivalence between LIG and MAR requires the completeness of the conditional model $f(y|x;\theta)$ for almost all $x \in \mathcal{X}$. If $X$ is fully observable, then Theorem 3.1 still holds. However, if $X$ may also be missing, the equivalence of MAR and LIG requires more strongly that the joint density of $(Y,X)$ belong to a complete parameter family. This situation has already been included in the above discussion, since some components of $Y$ can be treated as covariates. However, caution must be taken for the model parameterization, as in general the parameterization for the joint distribution of $(Y,X)$ is distinct from that for the conditional distribution of $Y$ on $X$.

REMARK 5.4.   A special case occurs where $Y$ is a scalar random variable and no covariates are involved in the model. Here $r$ is just 0 or 1, and MAR requires that $f(0|y;\psi)$ is independent of $y$. However, then $f(1|y;\psi) = 1 - f(0|y;\psi)$ must be independent of $y$ too, and so $Y$ and $R$ are statistically independent in the usual sense. This is the missing completely at random (MCAR) condition [Rubin (1976)]. So in this special case the conclusion of Theorem 3.1 is that for complete families

$$\text{LIG} \Leftrightarrow \text{MAR} \Leftrightarrow \text{MCAR}.$$

## APPENDIX

PROOF OF LEMMA 3.1. Suppose that $f_1(y^{(1)}; \phi_1(\theta))$ is not complete. Then there exists $t(y^{(1)}) \neq 0$ such that

$$\int t(y^{(1)}) f_1\{y^{(1)}; \phi_1(\theta)\} \, dy^{(1)} = 0 \qquad \forall \, \theta \in \Theta \text{ (or, equivalently, } \forall \, \phi_1 \in \Phi_1).$$

Then

$$\int t(y^{(1)}) f(y; \theta) \, dy$$

$$= \int t(y^{(1)}) f_1(y^{(1)}; \phi_1) \int f_{2|1}\{y^{(2)}|y^{(1)}, \phi_2(\theta; y^{(1)})\} \, dy^{(2)} \, dy^{(1)}$$

$$= \int t(y^{(1)}) f_1\{y^{(1)}; \phi_1(\theta)\} \, dy^{(1)}$$

$$= 0 \qquad \forall \, \theta \in \Theta,$$

contradicting the completeness of $f(y; \theta)$. Hence (a) is established.

Now suppose that (b) does not hold. Then there exists some $\mathcal{A} \subset \mathcal{Y}^{(1)}$ with nonzero probability under the marginal distribution of $Y^{(1)}$, such that for any $y^{(1)} \in \mathcal{A}$ the family of conditional densities $f_{2|1}\{y^{(2)}|y^{(1)}, \phi_2(\theta; y^{(1)})\}$, $\theta \in \Theta$, is not complete. There must then exist some function $w(y^{(1)}, y^{(2)}) \neq 0$ defined for $y^{(1)} \in \mathcal{A}$ and $y^{(2)} \in \mathcal{Y}^{(2)}$ such that

$$\int w(y^{(1)}, y^{(2)}) f_{2|1}\{y^{(2)}|y^{(1)}, \phi_2(\theta; y^{(1)})\} \, dy^{(2)} = 0 \qquad \forall \, \theta \in \Theta.$$

Now define

$$t(y) = \begin{cases} w(y^{(1)}, y^{(2)}), & y^{(1)} \in \mathcal{A}, \\ 0, & \text{otherwise.} \end{cases}$$

Clearly $t(y) \neq 0$, but

$$\int t(y) f(y; \theta) \, dy$$

$$= \int_{\mathcal{A}} f(y^{(1)}; \phi_1(\theta)) \int_{\mathcal{Y}_2} w(y^{(1)}, y^{(2)}) f_{2|1}\{y^{(2)}|y^{(1)}, \phi_2(\theta, y^{(1)})\} \, dy^{(2)} \, dy^{(1)}$$

$$= 0 \qquad \forall \, \theta \in \Theta.$$

This again contradicts the completeness of $f(y; \theta)$. $\square$

**Acknowledgments.** We are grateful to the referees for their helpful comments on an earlier draft of this paper.

## REFERENCES

ARNOLD, B. C., CASTILLO, E. and SARABIA, J. M. (1999). *Conditional Specification of Statistical Models.* Springer, New York. MR1716531

HEITJAN, D. F. (1993). Ignorability and coarse data: Some biomedical examples. *Biometrics* **49** 1099–1109.

HEITJAN, D. F. (1994). Ignorability in general incomplete-data models. *Biometrika* **81** 701–708. MR1326420

HEITJAN, D. F. (1997). Ignorability, sufficiency and ancillarity. *J. Roy. Statist. Soc. Ser. B* **59** 375–381. MR1440587

HEITJAN, D. F. and RUBIN, D. B. (1991). Ignorability and coarse data. *Ann. Statist.* **19** 2244–2253. MR1135174

JACOBSEN, M. and KEIDING, N. (1995). Coarsening at random in general sample spaces and random censoring in continuous time. *Ann. Statist.* **23** 774–786. MR1345200

KENWARD, M. G. and MOLENBERGHS, G. (1998). Likelihood based frequentist inference when data are missing at random. *Statist. Sci.* **13** 236–247. MR1665713

LITTLE, R. J. A. (1994). A class of pattern-mixture models for normal incomplete data. *Biometrika* **81** 471–483. MR1311091

LITTLE, R. J. A. and RUBIN, D. B. (2002). *Statistical Analysis with Missing Data*, 2nd ed. Wiley, New York. MR1925014

RUBIN, D. B. (1976). Inference and missing data (with discussion). *Biometrika* **63** 581–592. MR455196

SCHAFER, J. L. (1997). *Analysis of Incomplete Multivariate Data.* Chapman and Hall, London. MR1692799

TANNER, M. A. (1993). *Tools for Statistical Inference*: *Methods for the Exploration of Posterior Distributions and Likelihood Functions*, 2nd ed. Springer, New York. MR1238943

ZACKS, S. (1971). *The Theory of Statistical Inference.* Wiley, New York. MR420923

MRC HSRC
DEPARTMENT OF SOCIAL MEDICINE
UNIVERSITY OF BRISTOL
CANYNGE HALL
WHITELADIES ROAD
BRISTOL BS8 2PR
UNITED KINGDOM
E-MAIL: guobing.lu@bristol.ac.uk

DEPARTMENT OF STATISTICS
UNIVERSITY OF WARWICK
COVENTRY CV4 7AL
UNITED KINGDOM
E-MAIL: jbc@stats.warwick.ac.uk